\documentstyle[12pt,epsf]{article}
\begin{document}
\hspace\parindent
\thispagestyle{empty}

\bigskip
\bigskip
\bigskip

\centerline{\LARGE \bf Implementation of the Prelle-Singer}
\centerline{\LARGE \bf Method for 1ODEs}

\bigskip
\bigskip
\bigskip

\begin{center}
{\large
L.G.S. Duarte, S.E.S. Duarte, 
L.A.C.P. da Mota 
and J.E.F. Skea \footnote{E-mails: lduarte@dft.if.uerj.br, lduarte@dft.if.uerj.br, 
   damota@dft.if.uerj.br and jimsk@dft.if.uerj.br}
}

\end{center}

\bigskip

\centerline{\it Universidade do Estado do Rio de Janeiro,}
\centerline{\it Instituto de F\'{\i}sica, Depto. de F\'{\i}sica Te\'orica,}
\centerline{\it 20550-013 Rio de Janeiro, RJ, Brazil}

\bigskip
\bigskip

\begin{abstract} 

A set of MapleV R5 software routines for solving first order ordinary
differential equations (1ODEs) is presented. The package implements the 
Prelle-Singer Method in its original form (see \cite{PS}) plus its extension to include
elementary functions (ELFs)\footnote{For the definition of elementary function
see \cite{singer1}.} on the integrating factor (see \cite{Man1,Man2,ManMac}). The
package also presents a theoretical extension to deal with all 1ODEs presenting
liouvillian solutions (LIS) \cite{nossoPS1theo}.

\end{abstract}

\bigskip
\centerline{ \underline{\hspace{6.5 cm}} }

\medskip

\centerline{ {\bf (Submitted to Computer Physics Communications)} }

\bigskip
\bigskip
\bigskip
\bigskip
\bigskip
\bigskip

\newpage
\bigskip
\hspace{1pc}
{\bf PROGRAM SUMMARY}
\bigskip

\begin{footnotesize}
\noindent
{\em Title of the software package:} PSolver.   \\[10pt]
{\em Catalogue number:} (supplied by Elsevier)                \\[10pt]
{\em Software obtainable from:} CPC Program Library, Queen's
University of Belfast, N. Ireland (see application form in this issue)
\\[10pt]
{\em Licensing provisions:} none  \\[10pt]
{\em Operating systems under which the program has been tested:}
Linux (RedHat 5.2 and Debian 2.0.34), Windows 95, Windows 98.
\\[10pt]
{\em Programming languages used:} Maple V Release  5.
\\[10pt]
{\em Memory required to execute with typical data:}  32 Megabytes. \\[10pt]
{\em No. of lines in distributed program, including On-Line Help,
etc.:} 1370.                                                   \\[10pt]
{\em Nature of mathematical problem}\\
Solving of first order differential equations via the Prelle-Singer Method. 
\\[10pt]
{\em Methods of solution}\\
The method of solution is partially the standard Prelle-Singer Method extended to include
the case where the 1ODE presents elementary functions. We have also implemented
an extension of our own to deal with 1ODEs which solutions are liouvillian.
\\[10pt]
{\em Restrictions concerning the complexity of the problem}\\
If the integrating factor for the 1ODE in charge has factors of high degree in
the dependent and independent variables and in the ELFs appearing on the 1ODE
then the package may expend much computational time to catch it.
\\[10pt]
{\em Typical running time}\\
This depends strongly on the 1ODE. 
\\[10pt]
{\em Unusual features of the program}\\
Apart from being a solver of 1ODEs, implementing the Prelle -Singer approach, 
our package is a research tool allowing the user to follow each and every step
of the procedure, for example, which are the eigen-polynomials associated with 
the 1ODE (see section \ref{package}). Besides that, our package proved to be successful 
in solving many classes of 1ODEs that were missed by some of the best 
solvers available. Last but not least, our package implements a theoretical 
extension (for details, see section \ref{EPS} and our paper\cite{nossoPS1theo}) 
to the original Prelle-Singer approach that enhances its scope 
allowing for it to tackle all 1ODEs that present LIS.
\end{footnotesize}
\newpage
\hspace{1pc}
{\bf LONG WRITE-UP}

\section{Introduction}
\label{intro}

The fundamental position of differential equations (DEs) in scientific
progress has, over the last three centuries, led to a vigorous search
for methods to solve them. The overwhelming majority of these methods
are based on classification of the DE into types for which a method
of solution is known, which has resulted in a gamut of methods that
deal with specific classes of DEs. This scene changed somewhat at the
end of the 19th century when Sophus Lie developed a general method to
solve (or at least reduce the order of) ordinary differential equations
(ODEs) given their symmetry transformations ~\cite{step,bluman,olver}.  Lie's method
is very powerful and highly general, but first requires that we find
the symmetries of the differential equation, which may not be easy to do.
Search methods have been developed~\cite{nosso,nosso2} to extract
the symmetries of a given ODE, however these methods are heuristic and
cannot guarantee that, if symmetries exist, they will be found.

On the other hand in 1983 Prelle and Singer (PS) presented a deductive method
for solving 1ODEs that presents a solution in terms of ELFs (SELF) if such a
solution exists~\cite{PS}. The attractiveness of the PS method lies not only in
the fact that it is based on a totally different theoretical point of view but
also because, if the given 1ODE has a SELF, the method guarantees that this
solution will be found (though, in principle it can admittedly take an infinite
amount of time to do so).

Despite all the good points about the PS-method it presents a shortcoming: it
can not guarantee to solve 1ODEs which solutions are not SELF. 
In order to improve that situation, we have developed \cite{nossoPS1theo} an
extended PS-method. 

So, apart from implementing the PS-method (for the first time, in Maple, with
the complete approach and the additional bonus of the afore mentioned research
tools), we present and implement a theoretical extension to the scope of the method.

The paper is organized as follows: In section \ref{PS}, we present a short
introduction to the PS-approach. In the following section, we introduce the main
ideas of our extension to the procedure. Next, in section \ref{package},
we present the package. In sections \ref{examples} and \ref{perform} we, respectively, present some
examples and the performance analysis of the package. Finally, we end with the
conclusion.

\section{The Prelle-Singer Method}
\label{PS}

Despite its usefulness in solving 1ODEs, the Prelle-Singer procedure is
not very well known outside mathematical circles, and so we present
a brief overview of the main ideas of the procedure.

\subsection{Rational 1ODEs}
\label{rational}

Consider the class of 1ODEs which can be written as
\begin{equation}
\label{1ODE}
y' = {\frac{dy}{dx}} = {\frac{M(x,y)}{N(x,y)}}
\end{equation}
where $M(x,y)$ and $N(x,y)$ are polynomials with coefficients in the complex 
field $\it C$.

In~\cite{PS}, Prelle and Singer proved that, if an elementary first
integral of~(\ref{1ODE}) exists, it is possible to find an integrating
factor $R$ with $R^n~\in~\it C$ for some (possible non-integer) $n$, such that
\begin{equation}
{\frac{\partial (RN)}{\partial x}}+{\frac{\partial (RM)}{\partial y}} = 0.
\label{eq_int_factor}
\end{equation}

\noindent
The ODE can then be solved by quadrature.  From~(\ref{eq_int_factor})
we see that
\begin{equation}
  N\, {\frac{\partial R}{\partial x}}
+ R\, {\frac{\partial N}{\partial x}}
+ M\, {\frac{\partial R}{\partial y}}
+ R\, {\frac{\partial M}{\partial y}} 
= 0.
\label{eq_int_factor_aberta}
\end{equation}

\noindent
Thus
\begin{equation}
{\frac{D[R]}{R}} =  - \left( {\frac{\partial N}{\partial x}} + 
{\frac{\partial M}{\partial y}} \right),
\label{eq_PS}
\end{equation}
where
\begin{equation}
\label{eq_def_D}
D \equiv N {\frac{\partial }{\partial x}}
        + M {\frac{\partial }{\partial y}}.
\end{equation}

\noindent
Now let 
\begin{equation}
\label{R}
R = \prod_i f^{n_i}_i
\end{equation}
where $f_i$ are irreducible polynomials 
and $n_i$ are non-zero rationals. From (\ref{eq_def_D}), we have
\begin{eqnarray}
\label{ratio}
{\frac{D[R]}{R}} & = & {\frac{D[\prod_{i} f^{n_i}_i]}{\prod_i f^{n_k}_k}} =
                    {\frac{\sum_i f^{n_i-1}_i n_i D[f_i] \prod_{j \ne i}
                     f_j^{n_j}}{\prod_k f^{n_k}_k}} \nonumber \\[3mm]   
                 & = &  \sum_i {\frac{f^{n_i-1}_i n_i D[f_i]}{f_i^{n_i}}} =
                    \sum_i {\frac{n_iD[f_i]}{f_i}}.
\end{eqnarray}

From~(\ref{eq_PS}), plus the fact that $M$ and $N$ are polynomials, 
we conclude that ${D[R]}/{R}$ is a polynomial. Therefore,
from~(\ref{ratio}), we see that $f_i | D[f_i]$ (i.e., $f_i$ is a 
divisor of $D[f_i]$).

We now have a criterion for choosing the possible $f_i$ (build all
the possible divisors of $D[f_i]$) and, by using~(\ref{eq_PS})
and~(\ref{ratio}), we have
\begin{equation}
\label{eq_ni}
\sum_i {\frac{n_iD[f_i]}{f_i}} = -\left( {\frac{\partial N}{\partial 
x}} + {\frac{\partial M}{\partial y}} \right).
\end{equation}

If we manage to solve~(\ref{eq_ni}) and thereby find $n_i$,
we know the integrating factor for the 1ODE and the problem is
reduced to a quadrature. Risch's algorithm~\cite{Risch} can then
be applied to this quadrature to determine whether a solution
exists in terms of elementary functions.

\subsection{Algebraic and Transcendental 1ODEs}
\label{transcendental}
\indent

Note that the procedure explained above assumed (in eq.\ref{1ODE}) that
$M(x,y)$ and $N(x,y)$ are polynomials in $x$ and $y$ over ${\bf C}[x,y]$. 

In order to allow the PS-method to solve differential equations with
transcendental or algebraic terms, Shtokhamer 
\cite{Shtokhamer} regarded the different transcendental and algebraic terms appearing in
$M(x,y)$ and $N(x,y)$ as new variables ($u_1,u_2,...,u_n$) to be used, in addition to the original
variables $x$ and $y$, in the construction of the polynomials $f_i$ described
above. 

Bearing that in mind, it is easy to notice that eq. \ref{eq_PS} will
become:
\begin{equation}
\label{eqU}
{\frac{D[R]}{R}} = - \left({\frac{\partial}{\partial x}} + 
       \sum_{i=1}^{n} {\frac{\partial u_i}{\partial x}} 
          {\frac{\partial}{\partial u_i}} \right)
          N(x,y,u_{1...n}) - \left({\frac{\partial}{\partial y}} + 
          \sum_{i=1}^{n} {\frac{\partial u_i}{\partial y}}
          {\frac{\partial}{\partial u_{i}}} \right) M(x,y,u_{1...n})
\end{equation}

A quick inspection of the equation above reveals that its right-hand-side will 
not generally be polynomial, as it would be necessary to apply the PS-method. 

The non-polynomial terms might have two origins:
\begin{itemize}
\item ${\frac{\partial u_i}{\partial x}}$ and ${\frac{\partial u_i}{\partial y}}$ 
may produce functions that are not in $\{u_1,...,u_n\}$.
\item Further than that, ${\frac{\partial u_i}{\partial x}}$ and 
${\frac{\partial u_i}{\partial y}}$ could be rational instead of polynomial.
\end{itemize}

We can address the first question by noticing that the $u_i$, where $i=1,...n$,
present on the 1ODE are, by hypothesis, elementary. Therefore, we can construct
a complete set $\{u_1,..,u_n,..u_{n+m} \}$ (which we will call the basis of
functions), by including the new functions (generated by differentiation as
explained above) until we exhaust the process, i.e., no new functions appear.

The second problem can be fixed by multiplying the whole eq.(\ref{eqU}) by all the
denominators generated by ${\frac{\partial u_i}{\partial x}}$ and 
${\frac{\partial u_i}{\partial y}}$.

To illustrate these mystical arguments let us use a simple example.
Suppose we have the following 1ODE:

\begin{equation}
\label{eqex}
y'={\frac{ln(x)+sin(x)}{y(x)}}.
\end{equation}
So, we have
\begin{equation}
\label{us}
u_1=ln(x), \,\,\, u_2=sin(x).
\end{equation}
Let us first check if the derivatives introduce new functions:
\begin{equation}
\label{dus}
{\frac{\partial u_1}{\partial x}} = {\frac{1}{x}}, \,\,\, 
{\frac{\partial u_2}{\partial x}} = cos(x), \,\,\, {\frac{\partial u_1}{\partial
y}} = 0, \,\,\, {\frac{\partial u_2}{\partial y}} = 0.
\end{equation}

\noindent
From eq.(\ref{dus}), one can see that a new function, $cos(x)$, was introduced.
So, we have
\begin{equation}
\label{us2}
u_1=ln(x), \,\,\, u_2=sin(x), \,\,\, u_3=cos(x).
\end{equation}

Now, we would have to, once again, verify if the derivation will introduce new
functions. In this case, it is easy to see that no more functions will appear,
since 
\begin{equation}
\label{complete}
{\frac{\partial \,cos(x)}{\partial x}} = -sin(x), \,\,\, ({\frac{\partial
u_3}{\partial x}} = -u_2), 
\end{equation}
thus completing the set.

Therefore, we have solved the first problem mentioned above. From eq.(\ref{dus})
we can see that the derivative of $ln(x)$ generates a denominator $(x)$. 
Then, to ensure that the right-hand-side of eq.(\ref{eqU}) is polynomial, we would  have
to multiply the whole equation by $x$ (i.e., by the denominator of 
${\frac{\partial u_1}{\partial x}}$).

So, we can summarize the approach described in this section by three steps:
\begin{itemize}
\item complete the set of functions appearing on the 1ODE.
\item define a new differential operator ${\cal D}$ as ${\cal D} \equiv (\prod denominators) D$
\item substitute 
$
{\frac{\partial N}{\partial x}} + {\frac{\partial M}{\partial y}} \nonumber
$
by
$
\left( {\frac{\partial N}{\partial 
x}} + {\frac{\partial M}{\partial y}} \right) (\prod denominators). \nonumber
$
obtaining
\begin{equation}
\label{eq_ni2}
\sum_i {\frac{n_i{\cal D}[f_i]}{f_i}} = -\left( {\frac{\partial N}{\partial 
x}} + {\frac{\partial M}{\partial y}} \right) (\prod denominators).
\end{equation}
\item Then the procedure follows as in the case of the rational 1ODEs.
\end{itemize}

In this section we presented an overview of the PS procedure. Remember that we
can only guarantee to solve 1ODEs which solutions are SELF. In the next section, we
will present a theoretical extension of the PS approach that will permit us to
tackle all 1ODEs with LIS.

\section{Extended PS procedure}
\label{EPS}

In this section, we are going to briefly introduce the ideas behind a
theoretical extension\cite{nossoPS1theo} that will enhance the scope of applicability of the
PS-procedure. We are going to use an example to make clear what we are talking about.

Consider the following example, taken from the book by Kamke\footnote{
A classical reference on the field, generally used as a testing ground.}
 \cite{kamke}(eq. k21):
\begin{equation}
\label{Kamke21}
{\frac {dy}{dx}}=y^{2}-y\sin(x)+\cos(x) ,
\end{equation}
which solution is
\begin{equation}
\label{solKamke21}
y=\sin(x)-{\frac {{e^{-\cos(x)}}}{{\it C}+\int \!{e^{-\cos(x)}}{
dx}}} .
\end{equation}

Note that the solution is not a SELF. So the PS-procedure can not guarantee to
get it. To understand why, let us then consider the integrating factor for this
1ODE:

\begin{equation}
\label{RR}
R = {\frac {{e^{-\cos(x)}}}{\left (y-\sin(x)\right )^{2}}}
\end{equation}

What can we learn from that? Examining $R$ we can see that it is written as a
product of ELFs. Indeed, we can claim the $R$ is of the form:

\begin{equation}
\label{R2}
R = \prod_i h^{n_i}_i
\end{equation}
where, in this case, $h_1 = y-\sin(x)$ and $h_2 = e^{-\cos(x)}$. Let us now
compare this with eq.(\ref{R}). What is the fundamental difference between them?
The point is that, in eq.(\ref{R2}), not all $h_i$ are polynomials over the algebraic
extension of eq.(\ref{Kamke21}). In this example, the guilty party is $h_2$, thus
explaining why the PS-procedure would not work.

Let us introduce our extended PS-procedure that will allow to deal with all
those cases. Some of what is about to come will involve the Lie
theory\cite{step,bluman,olver} for 1ODEs. The reader that is not familiar with
this should not despair. We are going to explain (in English) what every step
means.

Every 1ODE of the form
\begin{equation}
\label{1ODE2}
y' = {\frac{dy}{dx}} = {\frac{M(x,y)}{N(x,y)}}
\end{equation}
where $M(x,y)$ and $N(x,y)$ are polynomials of ELFs,  presents an 
infinite number of symmetries. We can guarantee that at
least one is of the form:
\begin{equation}
\label{sym}
\left[\,\xi = 0, \,\,\eta = \eta(x,y)\,\right]
\end{equation}

In this language, the integrating factor is:
\begin{equation}
\label{R3}
R = {\frac{1}{M \xi - N \eta}} \,\,\,\,\Rightarrow \,\,\,\,R = {\frac{-1}{N \eta(x,y)}}
\end{equation}

So, if the integrating factor $R$, for the 1ODE (\ref{1ODE2}), is of 
the form given by eq.(\ref{R2}), then $\eta(x,y)$ must also be of the form
\begin{equation}
\eta(x,y) = \prod_k h^{m_k}_k,
\end{equation}
since, by assumption, N is polynomial on the algebraic extension 
of the 1ODE (\ref{1ODE2}). So, for the 1ODE (\ref{Kamke21}), using eqs.(\ref{RR}) and (\ref{R3}), we
have:
\begin{equation}
\label{ETA}
\eta(x,y) = - e^{\cos(x)}\left (y-\sin(x)\right )^{2}
\end{equation}

Applying the Lie theory one can obtain the most general form for an 1ODE, with a
given symmetry. For the symmetry presented on eq.(\ref{sym}), this is:

\begin{equation}
\label{eqgeneral}
{\frac {dy}{dx}}=\eta(x,y)\int \!{\frac {\eta_x(x,y)}{\eta(x,y)^{2}}}{dy}+
\eta(x,y) F(x)
\end{equation}

Examining eq.(\ref{Kamke21}), we see that it does not contain $e^{\cos(x)}$. On
the other hand, by looking at eqs.(\ref{ETA}) and (\ref{eqgeneral}), it apparently
should contain it! How can $e^{\cos(x)}$ disappear from eq.(\ref{eqgeneral})?!

First of all, from the last term of eq.(\ref{eqgeneral}), we can see that any
multiplicative factor, present on $\eta(x,y)$, can be
eliminated via the arbitrariness of $F(x)$. To be sure that this
multiplicative factor does not ``survive'' in the other term (the one with the
integral) it can not be a general function of $x$ and $y$. Actually, one can see (from
the general structure of the term) that this factor has to be an exponential
\cite{nossoPS1theo}.

Summarizing the ideas presented on this section, we see that the PS-procedure (for
both the rational and algebraic-transcendental cases) is based on the fact that
some quantities are divisors of others (see (\ref{eq_ni}) and (\ref{eq_ni2})).
With the example we have presented, we hope to have left it plain to the reader
that one might have factors appearing on $R$ that are not present on the 1ODE
and, indeed, that those factors have to be exponentials where the arguments are
functions of $x$ and $y$. The question that presents itself is: what
are those functions? By inspecting eqs. (\ref{eq_ni}) and (\ref{eq_ni2}) we can see
that they have to be taken from the basis of functions of the corresponding 1ODE
(defined, as explained on section \ref{transcendental}, by the functions present
on $N$ and $M$ of the 1ODE).

We would like to point out that the example dealt with on this section (and the
reasoning following it) is just a part of the whole body of ideas and arguments
we present on \cite{nossoPS1theo} to conclude that there is a possible extension to the
PS-procedure. Nevertheless, it is a complete part in the sense that it should be
enough to ``convince'' the reader of the importance of this extension. In other
words, we hope to have showed that, without this, a whole class of 1ODEs would
elude the PS-procedure. We will, later on this paper, show some more examples
only solved by this PS-extended-procedure and comment on our computational
implementation of these ideas.

\section{The {\it PSolver} package}
\label{package}

\subsection*{\it Summary}

A brief review of the commands of the package is as
follows:\footnote{This subsection and the next one may contain some
information already presented in the previous sections; this is
necessary to produce a self-contained description of the package.}

\begin{itemize}

\item {\tt PSolve} solves 1ODEs, using the Prelle-Singer approach. It deals not
only with rational 1ODEs but also with the ones containing ELFs (algebraic and
transcendental functions).

\item {\tt PSolvexp} extends the applicability of the PS-procedure to all 1ODEs
with LIS.

\item {\tt Intfact}  returns an integrating factor for the 1ODE.

\item {\tt Dop} constructs the ${\cal D}$ operator associated to the 1ODE.

\item {\tt Basis} constructs a list of equations of transformation which name a
set of algebraic and transcendental functions (e.g., $[u_1=f_1(x,y), ...]$).
Those functions are, see sections \ref{transcendental} and \ref{Basis}, what we call the
basis of functions for the corresponding 1ODE.

\item {\tt dBasis} presents the list with the partial derivatives of the basis
of functions (e.g., $[\partial u_1/\partial x = -u_2,...]$).

\item {\tt Eigenpol} calculates all the polynomials $f_i$ that are divisors of
${\cal D}f_i$, i.e., that satisfy the equation ${\cal D}f_i=g_i f_i$.

\item {\tt EigenPval} calculates  the polynomials $g_i$ that are the
``eigenvalues'' for the equation ${\cal D}f_i=g_i f_i$.

\end{itemize}

\subsection*{\it Description}

A complete description of the {\it PSolver} package's commands is
found in the on-line help. Here we present the most relevant part of
that description.

\subsection{Command name: {\tt PSolve}}
\label{PSolve}

\noindent {\it Feature:} This command applies the PS-procedure to solve 1ODEs.

\bigskip

\noindent
{\it Calling sequence\footnote{In what follows, the {\it input}
can be recognized by the Maple prompt \verb->-.}:}
\begin{verbatim}
> PSolve(1ode,n,optional_parameters);
\end{verbatim}

\noindent
{\it Parameters:}

\noindent
\begin{tabular}{ll}
\verb-1ode-           & - a first order ordinary differential equation.
\\
\verb-n-              & - an integer which sets the highest degree the command 
will consider for finding the \\ 
& polynomials $f_i$.
\end{tabular}

\bigskip

\noindent
{\it Optional Parameters:}
\smallskip
\smallskip

\noindent
\begin{tabular}{ll}
\verb-NumberF=integer-
   & - Tells the program the number of $f_i$ it is going to\\
   & use to look for the integrating factor. \\
\\
\verb-NumberF=All-
   & - Tells the program to use all $f_i$ in searching \\
   & for the integrating factor. \\
\\
\verb-AllSolutions-
   & - The absence of this argument truncates the searching for all\\
   &  possible $f_i$  in a way that, for the majority of cases, still allows  \\
   & the command to find the integrating factor but with considerable \\
   & gain in time.\\
\\   
\verb-Extended-  &  - Indicates to the program that it should use the set of
routines\\
&   that look for the exponentials  in order to calculate the integrating \\
& factor.\\
\\
\verb-Guess-
   & - It activates an algorithm for fast determination of $f_i$\\
   &  candidates.
\end{tabular}

\noindent

{\it Synopsis:}
\smallskip
\smallskip

The {\tt PSolve} command is a part of the {\bf PSolver} package, which, to
brutally summarize things, is designed to solve 1ODEs with LIS. The command implements
the PS-procedure and an extension of our own (see section \ref{EPS})
 which can be divided into five steps:
\begin{itemize}
\item First, from the 1ODE, one can determine the basis of functions.
\item In second place, from the basis (mentioned above) one constructs the $\cal
D$ operator, see section\ref{transcendental}.
\item Following that, we can calculate the ``eigen-polynomials'' and
``eigenvalues'' (which are themselves polynomials) for the
operator $\cal D$.
\item From these results, we find (or not, if it does not exists) the
integrating factor for the 1ODE.
\item Finally, one can then use the integrating factor to find the desired solution
for the 1ODE.
\end{itemize}

\medskip
\noindent
{\it The arguments}
\smallskip
\smallskip

The first argument of {\tt PSolve} is the 1ODE. The second argument
is an integer that tells the command the highest degree it is going to use in the
search for the integrating factor. We would like to remind the reader that the
PS-procedure is a semi-decision one\cite{PS,Man1}, in other words, if the integrating factor
has not be found for this degree, it does not exists at this level. But, it could
as well exist for higher degree, we have to keep looking.

\bigskip

\noindent
{\it The optional arguments}
\label{PSolveoptional}
\smallskip
\smallskip

The optional arguments can be given alone or in conjunction and in any
order. In order to improve its practicability (speed combined with efficiency)
{\tt PSolve}, by default, uses seven $f_i$ to construct the integrating
factor. Why is that so? Well, this process of calculate the integrating factor
can be broadly divided into two steps: first one has to find all possible $f_i$,
then use them to look for the integrating factor. This process of looking for
the IF is ``computationally expensive'' and, furthermore, we have noticed that it
is seldom the case where all the $f_i$ are needed to do the job. We have also
noticed that, generally, using the ``simpler'' $f_i$ was enough. What do we mean
by ``simpler''? We call a simple $f_i$ the one which leads to a faster
calculation. Luckily for us, there is a direct relation between this desired
property and the form of the $f_i$. Actually, we have observed that the
``smaller'' the $f_i$ the faster is the computation time. So we introduced a
algorithm for ordering the possible $f_i$ by ``size''\footnote{This definition of
``size'' and the creation of the routines to measure it are a part of our
package.}. The default number seven was defined heuristically. Of course there will
be the odd case where seven won't be enough, so we let the door open for the
user to increase this number via the optional argument {\tt NumberF=integer}.

How about the other step? Can we improve the efficiency of our package in the
finding of all $f_i$? In fact, for the 1ODEs which involve radicals, the problem
of computational-time-expenditure can be even greater. For instance, in a
particular example, the program took several minutes to find some forty $f_i$.
It turned out that we developed a different algorithm that, although loosing
some possible solutions, still produces enough of them (for a great number of
cases) thus allowing the program to find the IF. For the remaining (few) cases,
the user may then apply the optional argument {\tt AllSolutions} which will tell
the program to use the first (slower) algorithm.

The importance of the option {\tt Extended} lies in the fact it
implements a body of new ideas making it possible for the user to extend the
PS-procedure to a greater range of 1ODEs, to cover all the ones with LIS 
(not only with SELFs). See section \ref{EPS}. 

Finally, the {\tt Guess} option implements an algorithm for fast pattern
recognition determining possible $f_i$ candidates of high degree, thus shortening
the computational time required by the PS-procedure. This algorithm will help in
some (difficult) cases but it complicates the general one and so it
is optional.

\noindent
\subsection{Command name: {\tt Intfact}}
\label{Intfact}

\noindent 
{\it Feature:} Calculates a integrating factor ($R$) for the 1ODE.
\smallskip
\smallskip
\smallskip

\noindent {\it Calling sequence:}
\begin{verbatim}
> Intfact(1ode,n,optional_parameter);
\end{verbatim}

\noindent
{\it Parameters:}
\smallskip
\smallskip

\noindent
The parameters {\tt 1ode} and {\tt n} have the same meaning as explained
above\footnote{In what follows, these parameters will always have this meaning.}.

\bigskip

\noindent
{\it Optional Parameter:}
\smallskip
\smallskip

\noindent
\begin{tabular}{ll}
\verb-NumberF=integer-
   & - Tells the program the number of $f_i$ it is going to\\
   & use to look for the integrating factor. \\
\\
\verb-NumberF=All-
   & - Tells the program to use all $f_i$ in searching \\
   & for the integrating factor. \\
\\
\verb-AllSolutions-
   & - The absence of this argument truncates the searching for all\\
   &  possible $f_i$  in a way that, for the majority of cases, still allows  \\
   & the command to find the integrating factor but with considerable \\
   & gain in time.\\
\\   
\verb-Extended-  &  - Indicates to the program that it should use the set of
routines\\
&   that look for the exponentials  in order to calculate the integrating \\
& factor.\\
\\
\verb-Guess-
   & - It activates an algorithm for fast determination of $f_i$\\
   &  candidates.
\end{tabular}

\medskip
\noindent
{\it Synopsis:}
\smallskip
\smallskip

Actually this command embodies the heart of the package and that is so for,
basically, two reasons: The brilliance of the PS-procedure lies exactly on the
calculation of the integrating factor. After that step, the finding of the
solution is no novelty. The second reason is that our extension to the
PS-procedure can be also found mostly on this command (via the optional
parameter mentioned above). Again, after we find the integrating factor, it is
straightforward to look for the solution to the 1ODE.

\noindent
\subsection{Command name: {\tt Basis}}
\label{Basis}

\noindent {\it Feature:} Determines the basis of functions for the 1ODE.

\smallskip
\smallskip
\smallskip

\noindent {\it Calling sequence:}
\begin{verbatim}
> Basis(1ode);
\end{verbatim}

\noindent
{\it Synopsis:}
\smallskip
\smallskip

\noindent 
With the somewhat more specialized user in mind, we implement a set
of commands to help the research with 1ODEs. The first
one is this command that generates, from the ELFs present on the 1ODE, a set of
functions that is complete under differentiation.
\noindent
\subsection{Command name: {\tt DBasis}}
\label{DBasis}

\noindent {\it Feature:} Returns a list with the partial derivatives of the
functions in the Basis.

\smallskip
\smallskip
\smallskip

\noindent {\it Calling sequence:}
\begin{verbatim}
> DBasis(1ode);
\end{verbatim}

\noindent
{\it Synopsis:}
\smallskip
\smallskip

This is another one of the commands to help in the research on the structure of
1ODEs.  It generates a table with all the first partial derivatives of the
Basis, written in terms of this Basis and of the independent and dependent
variables. This will be used to construct the ${\cal D}$ operator.

\noindent
\subsection{Command name: {\tt Dop}}
\label{Dop}

\noindent {\it Feature:} Returns the ${\cal D}$ operator (for the 1ODE) that is
a essential ingredient in applying the PS-procedure.

\smallskip
\smallskip
\smallskip

\noindent {\it Calling sequence:}
\begin{verbatim}
> Dop(1ode);
\end{verbatim}

\noindent
{\it Synopsis:}
\smallskip
\smallskip

This command returns the ${\cal D}$ operator associated with the 1ODE. The
operator is written in the basis of functions for the 1ODE (see
section \ref{examples}).
\noindent
\subsection{Command name: {\tt Eigenpol}}
\label{Eigenpol}

\noindent {\it Feature:} Returns the ``eigen-polynomials'' for the ${\cal D}$ operator.

\smallskip
\smallskip
\smallskip

\noindent {\it Calling sequence:}
\begin{verbatim}
> Eigenpol(1ode,n);
\end{verbatim}

\noindent
{\it Synopsis:}
\smallskip
\smallskip

This command finds all the polynomials $f_i$ (of maximum degree {\tt n})
satisfying the equation ${\cal D}f_i = g_i f_i$ and those will be the candidates
to build the integrating factor for the 1ODE.

\noindent
\subsection{Command name: {\tt EigenPval}}
\label{EigenPval}

\noindent {\it Feature:} Returns the ``eigenvalues'' for the ${\cal D}$ operator.

\smallskip
\smallskip
\smallskip

\noindent {\it Calling sequence:}
\begin{verbatim}
> EigenPval(1ode,n);
\end{verbatim}

\noindent
{\it Synopsis:}
\smallskip
\smallskip

This command finds the polynomials $g_i$ that play the role of eigenvalues for
the equation ${\cal D}f_i = g_i f_i$.

\section{Examples}
\label{examples}

In this section, we are going to exemplify some aspects of the routines in the
package and show that, besides being valid only on the basis of being an
implementation of the regular PS-procedure (and that this fact already greatly
enhances the solving of 1ODEs capabilities within the Maple environment), our
program proved to be a successful attempt of implementing news ideas, extending
even further the applicability of the PS-procedure.

In order to avoid boring the reader to death and to allow him/her to focus on the
relevant aspects, we are going to concentrate (try to) in a few 1ODEs that will
serve as examples for as many features as possible.

\subsection{Usage of the package commands}
\label{simples}

Let us use a simple 1ODE to show the commands in action. Consider the following 1ODE:
\begin{equation}
\label{simple}
{\frac {dy}{dx}}={\frac {y\left (\cos(x)+y{e^{-x}}+1\right )}
{\cos(x)}} 
\end{equation}

Usually, people only want to find the solution for the 1ODE. So, that would
simply require the input lines:
\begin{verbatim}
> eq := diff(y(x),x) = y(x)*(cos(x)+y(x)*exp(-x)+1)/cos(x):
> PSolve(eq,1);
\end{verbatim}
\begin{equation}
\ln (y(x)+{e^{x}})-\ln (\tan(1/2\,x)-1)+\ln (\tan(1/2\,x)+1)-
\ln (y(x))-{\it \_C1}=0
\end{equation}

Sometimes the user is interested on having a look at just the integrating factor, so
he/she could type:
\begin{verbatim}
> R := Intfact(eq,1);
\end{verbatim}
$$
\begin{array}{c}
\mbox {{\it For the 1ODE in the form}} \\
\\
\displaystyle{{\frac {d}{dx}}y(x)={\frac {y(x)\left (\cos(x){e^{x}}+y(x)+{e^{x}} 
\right )}{\cos(x){e^{x}}}}} \\
\\
\mbox {{\it the integrating factor will be}} \\
\\
\displaystyle{R = {\frac {1}{y\cos(x)\left (y+{e^{x}}\right )}}
}
\end{array}
$$

\noindent 
Please notice that the program informs the user the format in which it is
considering the 1ODE to be in, that is important since it may be the case where
the input 1ODE gets transformed by the simplification routines and is regarded
in a different shape (of course it is still the same 1ODE) and that will affect the
integrating factor.

The package can also provide information about the basis of functions (that can
be useful in analyzing the structure of the 1ODE).
\begin{verbatim}
> basis := Basis(eq);
\end{verbatim}
$$
basis = [u_{{1}}=\cos(x),u_{{2}}=\sin(x),u_{{3}}={e^{x}}]
$$
For the  user that is interested in the limits of applicability of the PS-
procedure or that is otherwise involved in  more pure research in the algebraic
properties of 1ODEs, that can prove to be a very valuable piece of information.
The same applies to the commands that follow.

For instance, the following command to be exemplified ($dBasis$) provides
essential information in order to apply the PS-procedure (please refer to the 
above definition of $basis$ and, for further detail and meaning, see section\ref{transcendental}):
\begin{verbatim}
> d_basis := dBasis(eq);
\end{verbatim}
$$
d\_basis =[{\frac {\partial }{\partial x}}u_{{1}}(x,y)=-u_{{2}},{\frac {
\partial }{\partial x}}u_{{2}}(x,y)=u_{{1}},{\frac {\partial }{
\partial x}}u_{{3}}(x,y)=u_{{3}},
$$
$$
{\frac {\partial }{\partial y}}u_{{1}
}(x,y)=0,{\frac {\partial }{\partial y}}u_{{2}}(x,y)=0,{\frac {
\partial }{\partial y}}u_{{3}}(x,y)=0]
$$

The above result is used by the following command which defines the ${\cal
D}$ operator.
\begin{verbatim}
> D_operator:= Dop(eq);
\end{verbatim}
$$
D\_operator = w \mapsto u_{{1}}u_{{3}}{\frac {\partial}{\partial x}}w+
y\left (u_{{1}}u_{{3}}+y+u_{{3}}\right ){\frac {\partial}{\partial y}}w-
u_{{2}}u_{{1}}u_{{3}}{\frac {\partial}{\partial u_{{1}}}}w+
$$
$$
{u_{{1}}
}^{2}u_{{3}}{\frac {\partial}{\partial u_{{2}}}}w+
{u_{{3}}}^{2}u_{{1}}{\frac {\partial}{\partial u_{{3}}}}w
$$
The user who decides to use this command has to bear in mind that
it returns a procedure. Ops, what does that mean? It defines a mapping to be
used later on by the user. To clarify what we mean, the best to do is to see it
in action in a simple case. Let us apply the recently defined $D\_operator$ 
to (for example) $u_1$. One would get:
\begin{verbatim}
> example := D_operator(u[1]);
\end{verbatim}
$$
example = -u_{{2}}u_{{1}}u_{{3}}
$$

Of course, one has to apply $D\_operator$ to quantities written on the basis,
i.e., as a function (for our example) of $x,y,u_1,u_2$ and $u_3$.

The next two commands are very specialized (for the user really into 1ODE
research). Please see sections \ref{PS} and \ref{package} for recalling the meaning of them. For our
particular example, we have:
\begin{verbatim}
> eigen_p := Eigenpol(eq,1);
\end{verbatim}
$$
eigen\_p = [y,{e^{x}},\cos(x),y+{e^{x}},\cos(x)-I\sin(x)]
$$

\noindent
\begin{verbatim}
> eigen_p_val := EigenPval(eq,1);
\end{verbatim}
$$
eigen\_p\_val = [\cos(x){e^{x}}+y+{e^{x}},\cos(x){e^{x}},-{e^{x}}\sin(x),y+\cos(x){e^{
x}},-I\cos(x){e^{x}}]
$$

Let us have a closer look at this results. The two lists above are consistent in
the sense that their order is related, i.e., the ``eigen-value'' associated with
$y$ is $\cos(x){e^{x}}+y+{e^{x}}$, the ``eigen-value'' associated with ${e^{x}}$
is $\cos(x){e^{x}}$ and so on.
In the language of the basis of functions, that would translate to, for example:
$$ {\cal D} [y] = y (u_1u_3+y+u_3)$$

As previously mentioned, we have used a single 1ODE to exemplify the usage of
the package commands, we hope to have clarify that usage and helped the reader
to understand further the PS-procedure. Now we are going to further exemplify
the package but with a different approach. We will show that our package
complements the solving capabilities of the existing solvers in Maple. We will divide that
into two main groups: the first talks about the fact that our implementation
solves 1ODEs, with SELfs, that are missed by the powerful solvers on the MAPLE
algebraic package. The second one, deals with our theoretical extension to the
PS-procedure (and our implementation of that) and tackles 1ODEs with LIS.

\subsection{1ODEs with SELFs}
\label{elfs}

The solver (in MAPLE) is the $dsolve$ command and it is consider by many the most
powerful solver of ordinary differential equations commercially available. It 
fails to solve the following two examples:

\begin{equation}
eq_1 = {\frac {dy}{dx}}={\frac {{e^{x}}}{y^{2}}}+9\,{
\frac {\left ({e^{x}}\right )^{2}}{y^{2}}}-6\,y{e
^{x}}+y^{4}
\end{equation}

\begin{equation}
eq_2 = {\frac {dy}{dx}}={\frac {\left (y^{2}\ln (
x)^{5}+4\,y\ln (x)^{3}+4\,\ln (x)+y^{2}\right )
y^{2}}{\left (y\ln (x)^{2}+2\right )^{2}x}}
\end{equation}

Our program finds the following integrating factors for the above 1ODEs:
$$
R_1 = -\left (-3\,{e^{x}}+{y}^{3}\right )^{-2}
$$
$$
R_2 = {\frac {1}{x{y}^{4}}}
$$
which, in turn, lead to the following solutions:
$$
sol_1 = -1+9\,x{e^{x}}-3\,x y^{3}+9\,{\it C_1}\,{e^{x}}-3\,{\it \_C1}\,y^{3}
=0
$$
$$
sol_2 = -y^{3}\ln (x)^{6}-6\,y^{2}\ln (x)^{4}-
12\,y\ln (x)^{2}-6\,\ln (x) y^{3}-8+6\,{\it C_2}\,y^{3}=0
$$

It is worth reminding the reader that the {\tt dsolve} command brings within it
many algorithms implementing the Lie Method for solving ordinary differential
equations \cite{nosso,nosso2}. But, for an infinity of cases (of which the two
equations above are an example), the symmetries can not be found by the
this heuristic approach. On the other hand, the PS-procedure, being of  a 
semi-decision nature, guarantees that, if the solutions exists it will be caught (it
is only a matter of computational power).

\subsection{1ODEs with LIS}
\label{lis}

In this subsection, we are going to exemplify the usage of the theoretical
extension to the PS-procedure we have explained on section \ref{EPS}, which allows
the procedure to be extended to deal with 1ODEs with LIS.

Here are two such examples:

\begin{equation}
eq_3 = {\frac {dy}{dx}}={\frac {y+1+{e^{y}}{x}^{4}}{{x}^{2}y}}
\end{equation}

\begin{equation}
eq_4 = {\frac {dy}{dx}}={\frac {-\cos(x)xy-\cos(x)x-\cos(x)+y+1+x{e^
{y}}}{1+xy}}
\end{equation}

Our program finds the following integrating factors for the above 1ODEs:
$$
R_3 = {\frac {1}{{x}^{2}{e^{y}}{e^{{x}^{-1}}}}}
$$
$$
R_4 = {\frac {1}{{e^{y}}{e^{\sin(x)}}}}
$$
which, in turn, lead to the following solutions:
$$
sol_3 = -{e^{y+x^{-1}}}{\it C_3}+1+y+{e^{y+x^{-1}}}\int \!e^{x+x^{-1}}dx
$$
$$
sol_4 = \left (-xy-1-x\right ){e^{-y-\sin(x)}}-\int \!x{e^{-\sin(x)}}{dx
}-{\it C_4}=0
$$

As previously mentioned, these equations would elude the regular PS-procedure.
Furthermore, these examples (again they exemplify an infinitude of such cases)
also escape solution via the {\tt dsolve} command. So, the implementation of the
extended PS-procedure \cite{nossoPS1theo}, besides enlarging even more the scope of soluble
1ODEs, now allows for the fact that all 1ODEs with LIS (which include the ones
with SELFs) can be treated in the context of a semi-decision procedure.

\section{Performance}
\label{perform}

\begin{table}[htb]
\begin{center}
\begin{tabular}{|c|c|}
\hline 
Quadrature, Exact & 1, 4, 7, 12, 17, 23, 26, 39, 57, 59, 61, 69, 70, 71, 73, 76, 89,\\
Separable & 90, 131, 150, 154, 200, 209, 223, 227, 229, 245, 248, 251, \\
& 263, 267, 270, 271, 273, 274, 285, 288, 289, 290, 298, 299,\\
&  300, 305, 307, 308, 309, 310, 311, 322, 330, 335, 336, 341,\\
&  348, 352, 353, 355, 356, 358, 360, 361, 366\\
\hline 
ODEs presenting &  10, 11, 16, 33, 34, 35, 49, 50, 51, 53, 54, 55, 56, 72, 74, 79, 80,\\
 Arbitrary Functions &  84, 85, 86, 87, 110, 126, 127, 128, 201, 202, 212, 219, 230, 250,\\
&268, 269, 330, 331, 365, 367\\
\hline 
ODEs with Special  & 13, 14, 24, 25, 30, 36, 37, 40, 43, 45, 47, 48, 63, 66, 82, 88, 95,\\
Functions in the  & 99, 100, 105, 107, 111, 121, 139, 144, 145, 146, 147, 157, 166,\\
Solutions & 168, 169, 176, 179, 203, 205, 206, 234, 237, 253, 265\\
\hline
ODEs with LIS  & \\
out of the scope&5,18, 20, 21, 22, 27, 28, 46, 83, 129, 133, 164, 235, 343, 351 \\
of the PS-method\footnote{fazer!!}& \\
\hline 
\end{tabular}
\caption{Classification of ODEs (in Kamke) that will 
not be considered on the procedure testing (the last row will be considered
on the testing of the Extended PS-procedure).}
\label{perftab}
\end{center}
\end{table}

In this section, in order to show the efficiency of the ideas hitherto exposed,
we are going to use the equations on the book by Kamke \cite{kamke}. Which
equations should we choose from the book? There are 367 1ODEs of first degree on
the book. In principle, this would be our testing arena. However, among those
equations there are, for example, some quadrature, some exact 1ODEs, etc. So, we
see no point in using those in our test, therefore we will exclude them (see
table \ref{perftab}). There is another important restriction one should consider
before embarking on the test, some 1ODEs are ``out-of-scope'', i.e., they can
not be solved by the PS-procedure. More specifically, we mean 1ODEs with
arbitrary functions, with solutions in terms of special functions and some 1ODEs
with LIS (see table \ref{perftab}). The latter case will be only excluded from
our test for the usual PS-procedure. Surely, we will include those in our test
for the Extended-PS (See section \ref{teps}).

\subsection{Testing the usual PS-procedure}
\label{tups}

From  table \ref{perftab} we see that our arena for testing the usual PS-procedure consists
of 215 1ODEs. Below, see table \ref{perftab2}, we present the results of the
application of our package to those 1ODEs. 

To complete our test we have to mention the 1ODEs not present on table
\ref{perftab2}. Those are solved using the optional arguments described on
section\ref{package}. The 1ODEs 62, 112, 113, 114, 115, 116, 333, 357 are caught
using the optional parameter AllSolutions (115, 333 need, besides that, option
NumberF=All). The 1ODEs 151, 182, 185, 211, 257, 327 are caught by the single
use of NumberF=All. The 1ODEs 38, 42, 44, 77, 78, 106 are solved using the
special algorithm that allows for fast solution of ``hard'' cases (described on
section \ref{package}). This algorithm is enabled by the use of the optional
parameter, {\tt Guess}. 

Finally, the 1ODEs 52, 81, 142, 173, 184, 186-189, 266, 292 exceed the allowed time in
our testing (500 seconds).

\begin{table}[htb]
\begin{center}
\begin{tabular}{|c|c|c|}
\hline 
 & Average Time  & ODEs \\
&of solution&\\ 
\hline
& & 2, 3, 6, 8, 9, 19, 29, 31, 58, 60, 64, 65, 67, 68, 75,\\
& &91, 92, 93, 94, 96, 97, 98, 101, 102, 103, 117, 118, 119, 120, \\
&& 122, 123, 124, 125, 130, 132, 134, 135, 136, 137, 138, 148, \\
&&149, 152, 153, 155, 156, 158, 159, 160, 161, 162, 165, 167,\\
&&171, 174, 175, 177, 178, 180, 183, 190, 191, 192, 193, 194,\\
&&197, 198, 199, 204, 207, 210, 213, 214, 215, 216, 218, 221,\\
N = 1 & 1.3 s & 222, 224, 225, 226, 228,231, 232, 233, 236, 238, 239, 240, \\
& &241, 242, 243, 244, 246, 247, 249, 252,254, 255, 256, \\
&& 258, 259, 260, 261, 262, 264, 272, 275, 276, 277, 279, 280, \\
& &281, 282, 283, 284, 286, 287, 291, 293, 294, 295, 296, 297, \\
& &301, 302, 303, 306, 313, 314, 315, 317, 318, 319, 321, \\ 
& &323, 324, 325, 326, 327, 328, 329, 332, 333, 334, 338, 339, \\
&&  342, 344, 345, 346, 347, 349, 350, 354, 357, 359, 362, 363, 364\\
\hline
N = 2  & 3.3 s & 15, 41, 104, 108, 109, 140, 141, 143, 163, 170, 195, 208,\\
&&  217, 220, 304, 312, 340 \\
\hline
N = 3  & 26.5 s & 172, 181, 278, 320\\
\hline
\end{tabular}
\caption{This table presents 
the average time it took our program to solve 1ODEs in the book by Kamke. 
The table is organized by the value of N, the degree of the $f_i$
used to construct the integrating factor.}
\label{perftab2}
\end{center}
\end{table}

\subsection{Testing the Extended PS-procedure}
\label{teps}

From table \ref{perftab}, we come to the conclusion that 15 1ODEs should be
considered for our testing procedure. These are 1ODEs that would escape
solution via the usual PS-procedure. The 1ODEs are numbered 5, 18, 20, 21, 22, 27, 28, 
46, 83, 129, 133, 164, 235, 343, 351. From those, all but two are solved by our
implementation. The 1ODEs 22, 46 exceed the allowed testing time of 500 seconds.

To conclude our section on the performance of the package, let us just comment
that, apart from the examples we gave on section \ref{examples}, which are
representatives of an infinitude of cases that are missed by the powerful
solver in Maple (dsolve), the 1ODEs numbered 42, 151, 152, 185, 257, 350, 351
are also solved by our package and missed by dsolve. Besides that, the
implementation of the theoretical extension allows us to substantially decrease
(in some cases) the time needed to solve the 1ODE. For instance, from the list 
of 1ODEs that exceed the 500 seconds time limit we have set,
the 1ODE 142 is solved by the  extended procedure in only 3.4 seconds.

\section{Conclusions}
\label{conclude}

In this paper, we presented an implementation of the PS-procedure. As mentioned before, 
this method has the advantage of being of a semi-decision nature, i.e., it 
guarantees that, if the solutions exists it will be caught (it is only a question  
of time). Furthermore, we have made available to the user all the intermediary
steps of the PS-procedure: the integrating factor, the $\cal D$ operator, the
basis of functions for the 1ODE etc (see section \ref{package}).

We have also implemented a theoretical extension we have developed\cite{nossoPS1theo}. This
enlarges the scope of PS-procedure to allow it to tackle all 1ODEs with LIS
(that contain the set of the 1ODEs with SELfs previously tackled by the usual
PS-procedure). As an additional bonus to all this, we have enhanced the
1ODE solving capabilities in Maple, since this implementation of our extension
enables the package to solve 1ODEs previously missed by the powerful solver in
Maple (dsolve).

As further development for this work, we are currently working on extensions to
the package to tackle 1ODEs containing arbitrary functions and on the
implementation of the Collins approach\cite{ManMac} (and references therein).

\end{document}